\documentclass[12pt]{amsart}
\usepackage{amscd,amsfonts,amsmath,amssymb,amsthm,eucal,latexsym}
\usepackage[dvips]{graphicx}
\usepackage[all]{xy}

\def\iT#1{{\itshape#1\/}}

\numberwithin{equation}{section}
\newtheorem{prop}{Proposition}[section]
\newtheorem{thm}[prop]{Theorem}
\newtheorem{lem}[prop]{Lemma}
\newtheorem{cor}[prop]{Corollary}
\theoremstyle{definition}
\newtheorem{definition}[prop]{Definition}
\newtheorem{ex}[prop]{Example}

\newtheorem{rem}[prop]{Remark}
\newtheorem{conjecture}[prop]{Conjecture}
\newtheorem{notation}[prop]{Notation}

\def\Def#1{Definition~\ref{#1}}

\def\Thm#1{Theorem~\ref{#1}}

\long\def\BTHM#1#2{\begin{thm}\LL{#1}#2\end{thm}}
\long\def\BDF#1#2{\begin{definition}\LL{#1}#2\end{definition}}

\long\def\BCOR#1#2{\begin{cor}\LL{#1}#2\end{cor}}
\long\def\BNOT#1#2{\begin{notation}\LL{#1}#2\end{notation}}

\long\def\BFIG#1#2#3{\begin{figure}[th]#3\caption{#2}\LL{#1}\end{figure}}

\def\Pr{\begin{proof}}
\def\rP{\end{proof}}

\def\boxx#1{\leavevmode\vbox{\hbox to 0pt{\hss\raise1.8ex\vbox
to 0pt{\vss\hrule\hbox{\vrule\kern.75pt\vbox{\kern.75pt\hbox{\tiny #1}\kern.75pt}\kern.75pt\vrule}\hrule}}}\relax}

\def\LL#1{\label{#1}\protect\boxx{#1}}
\def\LL#1{\label{#1}}

\def\iFF{{if and only if }}

\def\s-s{{self-similar}}
\def\fr{finitely ramified}

\def\fr{fi\-ni\-te\-ly ram\-i\-fied}
\def\frf{\fr\ fra\-ct\-al}

\def\frcs{\fr\ cell struct\-ure}

\def\Lp{Laplacian}
\def\Df{Dirichlet form}

\def\bjs{basilica Julia set}

\def\a{\alpha}
\def\A{\mathcal A}

\def\E{\mathcal E}
\def\F{\mathcal F}

\def\V{\mathcal V}

\def\<{\langle}
\def\>{\rangle}

\def\Trace{\text{\rm Trace}\vph_}

\def\Dom{\text{\rm Dom\,}}
\def\vph{^{\vphantom{{\ }_{Ap}}}}

\DeclareMathOperator{\dom}{dom}%
\DeclareMathOperator{\diam}{Diam}%

\begin{document}

\title[Laplacians  on the basilica Julia set]
{Laplacians  on the basilica Julia set}

\author{Luke G. Rogers}
\author{Alexander Teplyaev}
\thanks{Research supported in part by the NSF grant DMS-0505622}

\address{\noindent Department of Mathematics,
University of Connecticut, Storrs CT 06269-3009 USA}
\email{rogers@math.uconn.edu}\email{teplyaev@math.uconn.edu}
\date{\today}
\begin{abstract} We consider the basilica Julia set of the polynomial $P(z)=z^{2}-1$
and construct all possible resistance (Dirichlet) forms, and the corresponding Laplacians, for
which the topology in the effective resistance metric coincides with the usual topology. Then we
concentrate on two particular cases. One is a self-similar harmonic structure, for which the energy
renormalization factor is $2$, the spectral dimension is $\log9/\log6$, and we can compute all the
eigenvalues and eigenfunctions by a spectral decimation method. The other is graph-directed
self-similar under the map $z\mapsto P(z)$; it has energy renormalization factor $\sqrt2$ and
spectral dimension $4/3$, but the exact computation of the spectrum is difficult. The latter
Dirichlet form and Laplacian are in a sense conformally invariant on the basilica Julia set.
\tableofcontents
\end{abstract}




\maketitle

\section{Introduction}

In the rapidly developing theory of analysis on fractals, the principal examples are \fr\ \s-s
fractal sets that arise as fixed points of iterated function systems (IFS).  For example, the
recent book of Strichartz \cite{Strichartz-book} gives a detailed account of the rich structure
that has been developed for studying differential equations on the the well known Sierpinski Gasket
fractal, primarily by using the methods of Kigami (see \cite{Ki}).
Some generalizations exist to fractals generated by graph-directed IFS and certain random IFS
constructions \cite{Hambly,HN,HMT}, but it is desirable to extend the methods to other interesting
cases. Among the most important and rich collections of fractals are the Julia sets of complex
dynamical systems (see, for instance, \cite{CG,DH,Milnor}).  In this paper we construct Dirichlet forms
and Laplacians on the Julia set of the quadratic polynomial $P(z)=z^{2}-1$, which is often referred
to as the \bjs, see Figure~\ref{basilicaj}. The \bjs\ is particularly interesting because it is one
of the simplest examples of a Julia set with nontrivial topology, and analyzing it in detail shows
how to transfer the differential equation methods of \cite{Strichartz-book} to more general Julia
sets.

\BFIG{basilicaj}{The \bjs, the Julia set of $z^2-1$.}%
{\includegraphics{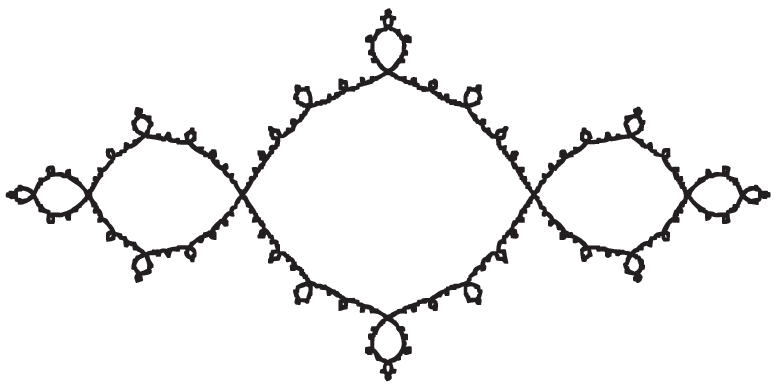}}

Another reason for our interest in the \bjs\ comes from its appearance as the limit set of the
so-called basilica self-similar group. This class of groups came to prominence because of their
relation  to finite automata and groups of intermediate growth, first discovered by Grigorchuk. The
reader can find extensive background on \s-s groups in the monograph of Nekrashevych
\cite{Nekrashevych}, some interesting calculations particularly relevant to the basilica group in
\cite{Kaimanovich}, and a review of the most recent developments in \cite{NekandTep}.

Since local regular Dirichlet forms and their \Lp s are in one-to-one correspondence, up to a
natural equivalence, with symmetric continuous diffusion processes and their generators, our
analysis allows the construction of diffusion processes on the \bjs. Random processes of this type
are interesting because they provide concrete examples of diffusions with nonstandard behavior,
such as sub-Gaussian transition probabilities estimates. For a detailed study of diffusions on some
\fr\ \s-s fractals see \cite{Barlow} and references therein.  The background on Dirichlet forms and
Markov processes can be found in \cite{FOT}.

Our construction starts with providing the \bjs\ with a finitely ramified cell structure,
(Definition~\ref{d-fract} in Section~\ref{frcssection}). According to \cite{T07cjm}, such a cell
structure makes it possible to use abstract results of Kigami \cite{Ki4} (see also \cite{Ki2008}) to construct local regular
\Df s that yield a topology equivalent to that induced from $\mathbb{C}$. In
Section~\ref{resformsection} we describe all such \Df s, and in Section~\ref{Lpsbjs} we describe
\Lp s corresponding to these forms and arbitrary  Radon measures. Among these \Lp s, some seem more
interesting than others. For example, there is a family of \Lp s with sufficient symmetry that
their graph approximations admit a spectral decimation like that in
\cite{Vibration,FukushimaShima,RammalToulouse,Shima}; this allows us to describe the eigenvalues and
eigenfunctions fairly explicitly in Section~\ref{specdecsection}. The energy renormalization factor
is $2$ and the spectral dimension is $\log9/\log6$ in this case.

In Section~\ref{sec-conf} we
describe the unique, up to a scalar multiple, \Df\ and \Lp\ that are conformally invariant under
the dynamical system. The latter \Lp\ does not have spectral decimation and we cannot determine its
eigenstructure, but its spectral dimension  can be computed to be $4/3$ using the renewal theorem
in \cite{HN}. The energy renormalization factor is $\sqrt2$ in the conformally invariant case. One
can find related group-theoretic computations and discussions in \cite{NekandTep}.

\section{A \frcs\ on the \bjs}\label{frcssection}

We will construct Dirichlet forms and Laplacians on the \bjs\ as limits of the corresponding
objects on a sequence of approximating graphs.  In order that we may later compare the natural
topology associated with the Dirichlet form with the induced topology from $\mathbb{C}$, we will
require some structure on these approximations. The ideas we need are from \cite{Ki4,T07cjm}, in
particular the following definition is almost identical to that in \cite{T07cjm}. The only change
is that we will not need the existence of harmonic coordinates, and therefore do not need to assume
that each $V_\alpha$ has at least two elements.

\BDF{d-fract}{A \iT{\fr\ set} $F$ is
a compact metric space with a \iT{cell structure}~$\F=\{F_\alpha\}_{\alpha\in\A}$ and
a \iT{boundary (\text{vertex}) structure} $\V=\{V_\alpha\}_{\alpha\in\A}$
such that
the following conditions hold.
\begin{itemize}
\item[(FRCS1)]
$\A$ is a countable index set;
 \item[(FRCS2)]
each $F_\alpha$
is a distinct compact connected subset of $F$;
 \item[(FRCS3)]
 each  $V_\alpha$ is a finite subset of
$F_\alpha$;
 \item[(FRCS4)]
if $F_\alpha=\bigcup_{j=1}^k F_{\alpha_j}$  then
$V_\alpha\subset\bigcup_{j=1}^k V_{\alpha_j}$;
 \item[(FRCS5)]
{there exists a filtration $\{\A_n\}_{n=0}^\infty$ such that
\begin{itemize}
\item[(i)]
$\A_n$ are finite subsets of $\A$,  $\A_{0}=\{0\}$, and
$F_{0}= F$;
\item[(ii)]
$\A_n\cap\A_m=\varnothing$ if $n\neq m$;
 \item[(iii)]
for any $\alpha\in\A_n$ there are
${\alpha_1},...,{\alpha_k}\in\A_{n+1}$
such that $F_\alpha=\bigcup_{j=1}^k F_{\alpha_j}$;
\end{itemize}}
 \item[(FRCS6)]
$F_{\alpha'}\bigcap F_{\alpha\vph}=V_{\alpha'}\bigcap V_{\alpha\vph}$ for any
two distinct $\alpha,\alpha'\in\A_n$;
\item[(FRCS7)]
for any strictly decreasing infinite cell sequence $F_{\alpha_1}\supsetneq
F_{\alpha_2}\supsetneq...$ there exists $x\in F$ such that
 $\bigcap_{n\geqslant1} F_{\alpha_n} =\{x\}$.
\end{itemize}
If these conditions are satisfied, then $$(F,\F,\V)=
(F,\{F_\alpha\}_{\alpha\in\A},\{V_\alpha\}_{\alpha\in\A})$$
is called a {\it \frcs.}}

\BNOT{not-Vn}{We denote $V_n=\bigcup_{\alpha\in\A_n}V_\alpha $.
Note that $V_n\subset V_{n+1}$ for all $n\geqslant0$ by \Def{d-fract}.
We say that $F_\alpha$ is an $n$-cell if $\alpha\in \A_n$.}

In this definition the vertex boundary $V_{0}$ of $F_{0}=F$ can be arbitrary, and in general may
have no relation with the topological structure of $F$.  However the cell structure is intimately
connected to the topology, as the following result shows.

\begin{prop}[\protect{\cite{T07cjm}}]\label{frcsandtopology}
The following are true of a \frcs.
\begin{enumerate}
\item For any $x\in F$ there is a strictly decreasing infinite sequence of cells
satisfying condition (FRCS7) of the definition. The diameter of cells in any such sequence tend to
zero.
\item The topological boundary of $F_\alpha$ is contained in $V_\alpha$ for any $\alpha\in\A$.
\item The set $V_*=\bigcup_{\alpha\in\A} V_\alpha$ is countably infinite, and $F$ is uncountable.
\item For any distinct $x,y\in F$ there is $n(x,y)$ such that if $m\geqslant n(x,y)$ then any $m$-cell can not
contain both $x$ and $y$.
\item  For any $x\in F$ and $n\geqslant 0$, let $U_n(x)$ denote the union of
all $n$-cells that contain $x$. Then the collection of open sets $\mathcal U=\{U_n(x)
^\circ\}_{x\in F,n\geqslant 0}$ is a fundamental sequence of neighborhoods. Here $B^\circ$ denotes
the topological interior of a set $B$. Moreover, for any $x\in F$ and open neighborhood $U$ of $x$
there exist $y\in V_*$ and $n$ such that $x\in U_n(x)\subset U_n(y)\subset U$. In particular, the
smaller collection of open sets $\mathcal U'=\{U_n(x) ^\circ\}_{x\in V_*,n\geqslant 0}$ is  a
countable fundamental sequence of neighborhoods.
\end{enumerate}
\end{prop}

In general a \frf\ may have many filtrations, and the Dirichlet forms, resistance forms and energy
measures we will discuss later are independent of the filtration. However it is natural in the
context of a self-similar set to consider a filtration that is adapted to the self-similarity.  We
now define a \frcs\ and a filtration, which have certain self-similarity properties, on the \bjs.


By definition, the 0-cell is the \bjs\ fractal, which we denote by $J$. Let us write
$a=\frac{1-\sqrt5}{2}$ for one of the fixed points of $z^{2}-1$. The interiors of four 1-cells are
obtained by removing the points $\pm a$; \BFIG{basilicaj2}{The  the Julia set of $z^2-1$ with the repulsive fixed points $a=\frac{1-\sqrt5}{2}$ and $b=\frac{1+\sqrt5}{2}$ circled.}%
{\begin{picture}(222.5, 106.5)(-111.25, 0)\put(-111.25,0){\includegraphics{basilicaj.eps}} \put(42,53.75){\setlength{\unitlength}{1pt}\circle{6}}
\put(-42,53.75){\setlength{\unitlength}{1pt}\circle{6}}
\put(110.5,53.75){\setlength{\unitlength}{1pt}\circle{6}}
\put(-110.5,53.75){\setlength{\unitlength}{1pt}\circle{6}}
\put(0,53.75){\setlength{\unitlength}{1pt}\circle*{2}}
\thicklines
\put(3,50 ){$0$}
\put(-37,50 ){$a$}
\put(37,50 ){\llap{$-a$}}
\put(117,50 ){$b$}
\put(-117,50 ){\llap{$-b$}}
\end{picture}}this disconnects the part of $J$ surrounding the basin
around $0$ into symmetric upper and lower pieces, and separates these from two symmetric arms, one
on the left and one on the right, see Figure~\ref{basilicaj2} (and also Figure~\ref{Gn}). The top and bottom cells we denote
$J_{(1)}$ and $J_{(2)}$ respectively, and the left and right cells we denote $J_{(3)}$ and
$J_{(4)}$ respectively. The cells $J_{(1)}$ and $J_{(2)}$ each have two boundary points, while
$J_{(3)}$ and $J_{(4)}$ each have one boundary point.  In the notation of Definition \ref{d-fract},
\begin{equation*}
    V_{(1)}=V_{(2)}=\{\pm a\},\quad V_{(3)}=\{-a\},\quad V_{(4)}=\{a\}
    \end{equation*}
and therefore the boundary set of the fractal is $V_{0}=\{-a,a\}$. Note that the other fixed point, $b=\frac{1+\sqrt5}{2}$, does not play any role in defining the cell structure.

For $n\ge1$ we set $\A_n=\{1,2,3,4\}\times\{1,2,3\}^{n-1}$.  To define the smaller cells, we
introduce the following definition. If a cell has two boundary points, it is called an
\emph{arc-type cell}. If a cell has one boundary point, it is called a \emph{loop-type cell}.

\begin{figure}
\begin{picture}(120, 60)(-60, 0) \thicklines \put(-60,-7 ){$v_{\a1}$} \put(60,-7 ){$v_{\a2}$}
\put(0,23 ){$v_{\a3}$} \setlength{\unitlength}{30pt} \qbezier(0,1)(1,1)(2,0)
\qbezier(0,1)(1,2)(0,2) \qbezier(0,1)(-1,1)(-2,0) \qbezier(0,1)(-1,2)(0,2)
\put(0,1){\setlength{\unitlength}{1pt}\circle*{5}}
\put(2,0){\setlength{\unitlength}{1pt}\circle*{5}}
\put(-2,0){\setlength{\unitlength}{1pt}\circle*{5}}
\end{picture}
\caption{An arc-type cell}
\label{decompositionofarctypecell}
\end{figure}
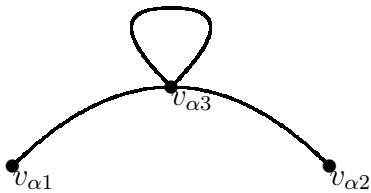

Each arc-type $n$-cell $J_\a$ is a union of three $n+1$-cells $J_{\a1}$,  $J_{\a2}$ and $J_{\a3}$;
$J_{\a1}$ and $J_{\a2}$ are arc-type cells connected at a middle point, while $J_{\a3}$ is a
loop-type cell attached at the same point (Figure~\ref{decompositionofarctypecell}).

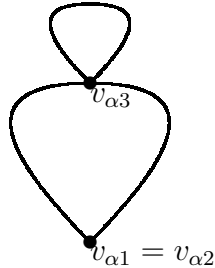
\begin{figure}
\begin{picture}(120, 90)(-60, 0) \thicklines \put(0,30 ){\put(0,23 ){$v_{\a3}$}
\setlength{\unitlength}{30pt} \qbezier(0,1)(1,2)(0,2) \qbezier(0,1)(-1,2)(0,2)
\put(0,1){\setlength{\unitlength}{1pt}\circle*{5}}} \put(0,-7 ){$v_{\a1}=v_{\a2}$} \put(0,-60
){\setlength{\unitlength}{60pt} \qbezier(0,1)(1,2)(0,2) \qbezier(0,1)(-1,2)(0,2)
\put(0,1){\setlength{\unitlength}{1pt}\circle*{5}}}
\end{picture}
\caption{A loop-type cell}\label{decompositionoflooptypecell}
\end{figure}

Each loop-type $n$-cell $J_\a$ is a union of three $n+1$-cells, $J_{\a1}$,  $J_{\a2}$ and
$J_{\a3}$; $J_{\a1}$ and $J_{\a2}$ are arc-type cells connected at two points, one of which is the
unique boundary point $v_{\a}\in V_{(\alpha)}$, while the other is the boundary point of the
loop-type cell $J_{\a3}$ (Figure \ref{decompositionoflooptypecell}).

The existence of this decomposition is a consequence of known results on the topology of quadratic
Julia sets.  In essence we have used the fact that the filled Julia set is the closure of the union
of countably many closed topological discs, and that the intersections of these discs are points
that are dense in $J$ and pre-periodic for the dynamics.  The Julia set itself consists of the
closure of the union  of the boundaries of these topological discs. This structure occurs for the
Julia set of every quadratic polynomial $z^{2}+c$ for which $c$ is in the interior of a hyperbolic
component of the Mandelbrot set or is the intersection point of two hyperbolic components, so in
particular for the \bjs\ because $c=-1$ lies within the period $2$ component. Details may be found
in \cite{CG,DH,Milnor}.  These general results imply that a FRCS may be obtained for all quadratic
Julia sets with suitable $c$ values in the manner similar to that described above, however the
\bjs\ is a sufficiently simple case that the reader may prefer to verify directly that the
existence (see, for instance, \cite[Theorem 18.11]{Milnor}) of internal and external rays landing
at $a$ implies that deletion of $\pm a $ decomposes $J$ into the four components $J_{(i)}$,
$i=1,2,3,4$, while the remainder of the decomposition follows by examining the inverse images of
these sets under the dynamics.

\begin{definition}\label{def-Gn}
The basilica  self-similar sequence of graphs  $G_n$ have vertices $V_{n}$ as previously described.
There is one edge for each pair of vertices joined by an arc-type $n$-cell, as well as one loop at
each vertex at which there is a loop-type $n$-cell. The result is shown in Figure~\ref{Gn}, and we
emphasize that it is highly dependent on our choice of filtration.
\BFIG{Gn}{Basilica self-similar sequence of graphs: graphs $G_0$ and $G_1$. }%
{\begin{picture}(80, 40)(-40, -20) \thicklines \setlength{\unitlength}{20pt}
\qbezier(0,1)(1,1)(2,0) \qbezier(2,0)(4,2)(4,0) \qbezier(0,-1)(1,-1)(2,0) \qbezier(2,0)(4,-2)(4,0)
\qbezier(-0,1)(-1,1)(-2,0) \qbezier(-2,0)(-4,2)(-4,0) \qbezier(-0,-1)(-1,-1)(-2,0)
\qbezier(-2,0)(-4,-2)(-4,0) \put(2,0){\setlength{\unitlength}{1pt}\circle*{5}}
\put(-2,0){\setlength{\unitlength}{1pt}\circle*{5}}
\put(-2.1,.3){$a$}
\put(1.6,.3){$-a$}
\end{picture}

\

\begin{picture}(80, 80)(-40, -40)
\thicklines \setlength{\unitlength}{20pt} \qbezier(0,1)(1,1)(2,0) \qbezier(0,1)(1,2)(0,2)
\qbezier(2,0)(4,2)(4,0) \qbezier(4,0)(5,1)(5,0)
\qbezier(0,-1)(1,-1)(2,0) \qbezier(0,-1)(1,-2)(0,-2) \qbezier(2,0)(4,-2)(4,0)
\qbezier(4,0)(5,-1)(5,0)
\qbezier(0,1)(-1,1)(-2,0) \qbezier(0,1)(-1,2)(0,2) \qbezier(-2,0)(-4,2)(-4,0)
\qbezier(-4,0)(-5,1)(-5,0)
\qbezier(0,-1)(-1,-1)(-2,0) \qbezier(0,-1)(-1,-2)(0,-2) \qbezier(-2,0)(-4,-2)(-4,0)
\qbezier(-4,0)(-5,-1)(-5,0)
\put(0,1){\setlength{\unitlength}{1pt}\circle*{5}}
\put(0,-1){\setlength{\unitlength}{1pt}\circle*{5}}
\put(2,0){\setlength{\unitlength}{1pt}\circle*{5}}
\put(-2,0){\setlength{\unitlength}{1pt}\circle*{5}}
\put(4,0){\setlength{\unitlength}{1pt}\circle*{5}}
\put(-4,0){\setlength{\unitlength}{1pt}\circle*{5}}
\put(-2.1,.3){$a$}
\put(1.6,.3){$-a$}
\end{picture}}\end{definition}

The sequence of graphs in Figure~\ref{Gn} is well adapted to the construction of the Kigami
resistance forms, and hence the Dirichlet forms, on $J$. For this reason it plays a prominent role
in Section~\ref{resformsection}.  In Section~\ref{specdecsection} a spectral decimation method for
this sequence of graphs is used to obtain a full description of the corresponding Laplacian.

It should be noted, however, that this is not the only sequence of graphs that we will consider. A
sequence that is arguably more natural, is the conformally invariant graph-directed sequence of
graphs for the \bjs, shown in Figure~\ref{d-c.i.g.}.

\BFIG{d-c.i.g.}{Basilica  conformally invariant graph-directed sequence of graphs. }%
{\begin{picture}(80, 40)(-40, -20) \thicklines \setlength{\unitlength}{20pt}
\qbezier(-2,0)(-4,-2)(-4,0)
\qbezier(-2,0)(-4,2)(-4,0)
\qbezier(-2,0)(4,-2)(4,0)
\qbezier(-2,0)(4,2)(4,0)
\put(-2,0){\setlength{\unitlength}{1pt}\circle*{5}}
\put(-2,.2){$a$}
\put(-3.5,-.2){$A$}
\put(2,-.2){$B$}
\end{picture}

\

\begin{picture}(80, 50)(-40, -20) \thicklines \setlength{\unitlength}{20pt}
\qbezier(0,1)(1,1)(2,0) \qbezier(2,0)(4,2)(4,0) \qbezier(0,-1)(1,-1)(2,0) \qbezier(2,0)(4,-2)(4,0)
\qbezier(-0,1)(-1,1)(-2,0) \qbezier(-2,0)(-4,2)(-4,0) \qbezier(-0,-1)(-1,-1)(-2,0)
\qbezier(-2,0)(-4,-2)(-4,0) \put(2,0){\setlength{\unitlength}{1pt}\circle*{5}}
\put(-2,0){\setlength{\unitlength}{1pt}\circle*{5}}
\put(-2.1,.3){$a$}
\put(1.6,.3){$-a$}
\put(-3.5,-.2){$B$}
\put(3.2,-.2){$B$}
\put(-.5,-.9){$A$}
\put(.5,1){$A$}
\end{picture}

\

\begin{picture}(80, 50)(-40, -20)
\thicklines \setlength{\unitlength}{20pt} \qbezier(0,1)(1,1)(2,0)
\qbezier(2,0)(4,2)(4,0) \qbezier(4,0)(5,1)(5,0)
\qbezier(0,-1)(1,-1)(2,0)  \qbezier(2,0)(4,-2)(4,0)
\qbezier(4,0)(5,-1)(5,0)
\qbezier(0,1)(-1,1)(-2,0)  \qbezier(-2,0)(-4,2)(-4,0)
\qbezier(-4,0)(-5,1)(-5,0)
\qbezier(0,-1)(-1,-1)(-2,0)  \qbezier(-2,0)(-4,-2)(-4,0)
\qbezier(-4,0)(-5,-1)(-5,0)
\put(2,0){\setlength{\unitlength}{1pt}\circle*{5}}
\put(-2,0){\setlength{\unitlength}{1pt}\circle*{5}}
\put(4,0){\setlength{\unitlength}{1pt}\circle*{5}}
\put(-4,0){\setlength{\unitlength}{1pt}\circle*{5}}
\put(-2.1,.3){$a$}
\put(1.6,.3){$-a$}
\put(-4.8,-.2){$B$}
\put(4.5,-.2){$B$}
\put(-.5,-.9){$B$}
\put(.5,1){$B$}
\put(-3.1,1){$A$}
\put(-2.7,-.9){$A$}
\put(3,1){$A$}
\put(2.2,-.9){$A$}
\end{picture}}

\BFIG{f-sub}{The substitution scheme for the basilica conformally invariant graph-directed sequence
of graphs. }{\begin{picture}(60, 70)(0, -10) \setlength{\unitlength}{10pt} \put(2.75,5 ){$A$}
\put(2.75,1 ){$B$} \put(3,4 ){\vector(0,-1){1.5}} \thicklines
\qbezier(0,5)(3,4)(6,5)\qbezier(0,1)(3,0)(6,1)
\end{picture}
\qquad\qquad
\begin{picture}(60, 70)(0, -10)
\setlength{\unitlength}{10pt} \put(2.75,5 ){$B$} \put(1.3,1 ){$A$}\put(4.3,1 ){$A$}\put(2.65,-.5
){$B$} \put(3,4 ){\vector(0,-1){1.5}} \thicklines
\qbezier(0,5)(3,4)(6,5)\qbezier(0,1)(1.5,0.25)(3,1) \qbezier(3,1)(4.5,0.25)(6,1)
\qbezier(3,1)(1,-1)(3,-1) \qbezier(3,1)(5,-1)(3,-1)
\end{picture} }

These graphs will be considered in Section~\ref{sec-conf}, where detailed definitions can be found.
Their construction is related to group-theoretic results \cite[and references
therein]{Nekrashevych,NekandTep}, and in particular to the substitution scheme in
Figure~\ref{f-sub}. The cell structure and the filtration could be defined starting with the single
point boundary set $\{a\}$, and then taking the inverse images $P^{-n}\{a\}$ of this point under
the polynomial $P(z)$, which is of course different from  $V_n$  in
Definition~\ref{def-Gn}. More precisely, for any $n$ and $k$ we have $V_n\neq P^{-k}\{a\}$, even though $V_*=\bigcup_{n\geq0}V_n=\bigcup_{n\geq0} P^{-n}\{a\}$.

\section{Kigami's resistance  forms  forms on the \bjs\ and the local resistance
metric}\label{resformsection}

One way of constructing Dirichlet forms on a fractal is to take limits of resistance forms on an
approximating sequence of graphs.  We recall the definition from \cite{Ki}, as well as the
principal results we will require.

\BDF{d-resistance form}{A pair $(\E,\Dom\E)$ is called a resistance form on a countable set $V_*$
if it satisfies the following conditions.
\begin{itemize}
\item[(RF1)] $\Dom\E$ is a linear subspace of $\ell(V_*)$ containing constants,
 $\E$ is a nonnegative symmetric quadratic form on $\Dom\E$, and $\E(u,u)=0$ if and only if $u$ is constant on $V_*$.
\item[(RF2)] Let $\sim$ be the equivalence relation on $\Dom\E$ defined by $u\sim v$ \iFF $u-v$ is constant on $V_*$. Then $(\E/\mbox{\hskip-.4em}\sim,\Dom\E)$ is a Hilbert space.
\item[(RF3)] For any finite subset $V\subset V_*$ and for any
$v\in\ell(V)$ there exists $u\in\Dom\E$ such that $u\big|_V=v$.
\item[(RF4)] For any $x,y\in V_*$ the resistance between $x$ and $y$ is defined to be
$$
R(x,y)=\sup\left\{ \frac{\big(u(x)-u(y)\big)^2}{\E(u,u)}:u\in\Dom\E,\E(u,u)>0\right\}<\infty.
$$
\item[(RF5)] For any $u\in\Dom\E$ we have the $\E(\bar u,\bar u)\leqslant\E(u,u)$, where
$$
\bar u(x)=\left\{
\begin{aligned}
&1&& \text{if \,} u(x)\geqslant1,\\
&u(x)&& \text{if \,} 0<u(x)<1,\\
&0&& \text{if \,} u(x)\leqslant1.
\end{aligned}\right.
$$
Property~(RF5) is called the Markov property.
\end{itemize}}

\begin{prop}[Kigami, \protect{\cite{Ki4}}]\label{basicfactsaboutresistmetric} Resistance forms have the following properties.
\begin{enumerate}
\item The effective resistance $R$ is a metric on $V_*$.  Any function in
$\Dom\E$ is $R$-continuous; in particular, if $\Omega$ is the $R$-completion of $V_*$ then any
$u\in\Dom\E$ has a unique $R$-continuous extension to $\Omega$.
\item For any finite subset $U\subset V_*$, a finite dimensional Dirichlet form $\E_U$ on $U$ may be defined
by $$\E_U(f,f)=\inf\{\E(g,g):g\in\Dom\E, g\big|_U=f\}.$$  There is a unique $g$ at which the
infimum is achieved. The form $\E_U$ is called the trace of $\E$ on $U$, and may be written
$\E_U=\Trace U(\E)$.  If $U_1\subset U_2$ then $\E_{U_1}=\Trace {U_1}(\E_{U_2})$.
\end{enumerate}
\end{prop}

Our description of the Dirichlet forms on the \bjs\ relies on the following theorems.

\begin{thm}[Kigami, \protect{\cite{Ki4}}]\label{Dirformdeterminedbytraces}
Suppose that $V_n$ are finite subsets of $V_*$ and that $\bigcup_{n=0}^\infty V_n$ is $R$-dense in
$V_*$. Then $$\E(f,f)=\lim_{n\to\infty}\E_{{V_n}}(f,f)$$ for any $f\in\Dom\E$, where the limit is
actually non-decreasing. Is particular, $\E$ is uniquely defined by the sequence of its finite
dimensional traces $\E_{V_n}$ on ${V_n}$.
\end{thm}

\begin{thm}[Kigami, \protect{\cite{Ki4}}]\label{Dirformconstructedfromtraces}
Suppose that  $V_n$ are finite sets, for each $n$ there is a resistance form $\E_{V_n}$ on ${V_n}$,
and this sequence of finite dimensional forms is compatible in the sense that each $\E_{V_n}$ is
the trace of $\E_{V_{n+1}}$ on ${V_n}$, where $n=0,1,2,...$. Then there exists a resistance form
$\E$ on $V_*=\bigcup_{n=0}^\infty V_n$ such that
$$\E(f,f)=\lim_{n\to\infty}\E_{{V_n}}(f,f)$$ for any $f\in\Dom\E$, and the limit is actually
non-decreasing.
\end{thm}

For convenience we will write $\E_n(f,f)=\E_{V_n}(f,f)$.  A function is called harmonic if it
minimizes the energy for the given set of boundary values, so a harmonic function is uniquely
defined by its restriction to $V_0$. It is shown in \cite{Ki4} that any function $h_0$ on  $V_0$
has a unique continuation to a harmonic function $h$, and $\E(h,h)=\E_n(h,h)$ for all $n$. This
latter is also a sufficient condition: if $g\in\Dom\E$ then $\E_0(g,g)\leqslant\E(g,g)$ with
equality precisely when $g$ is harmonic.

For any function $f$ on $V_{n}$ there is a unique energy minimizer $h$ among those functions equal
to $f$ on $V_{n}$. Such energy minimizers are called $n$-harmonic functions. As with harmonic
functions, for any function $g\in\Dom\E$ we have $\E_n(g,g)\leqslant\E(g,g)$, and $h$ is
$n$-harmonic if and only if $\E_n(h,h)=\E(h,h)$.

It is  proved in \cite{T07cjm} that if all $n$-harmonic functions are continuous in the topology of
$F$ then any $F$-continuous function is $R$-continuous and any $R$-Cauchy sequence converges in the
topology of $F$. In such a case there is also a continuous injection $\theta:\Omega\to F$ which is
the identity on $V_*$, so we can identify $\Omega$ with the the $R$-closure of $V_*$ in $F$. In a
sense, $\Omega$ is the set where the Dirichlet form $\E$ ``lives".

\begin{thm}[\protect{\cite{Ki4,T07cjm}}]
Suppose that all $n$-harmonic functions are continuous. Then  $\E$ is a local regular \Df\ on
$L^{2}(\Omega,\gamma)$, where $\gamma$ is any finite Borel measure on $(F, R)$ with the property
that all nonempty open sets have positive measure.
\end{thm}

\begin{proof}
The regularity of $\E$ follows from \cite[Theorem 8.10]{Ki4}, and its locality from  \cite[Theorem
3]{T07cjm}. Note that, according to \cite[Theorem 8.10]{Ki4}, in general for  resistance forms one
can consider $\sigma$-finite  Radon measures $\gamma$. However for a compact set a Radon measure
must be finite.
\end{proof}

The trace of $\E$ to the finite set $V_{n}$ may be written in the form
\begin{equation}\label{formonfiniteset}
    \E_{n}(f,f) = \sum_{\a\in \A_{n}} r^{-1}_{\a} \bigl( f(v_{\a1})-f(v_{\a2})\bigr)^{2},
    \end{equation}
from which we define the resistance across $J_{\a}$ to be the value $r_{\a}$.  Note that it is not
the same as $R(v_{\a1},v_{\a2})$.

The values $r_\a$ may be used to define a geodesic metric that is comparable to the resistance
metric. A {\em path} from $x$ to $y$ in $\Omega$ consists of a doubly infinite sequence of vertices
$\{v_{\a_{j}}\}_{-\infty}^{\infty}$ and arc-type cells $J_{j}$ connecting $v_{\a_{j}}$ to
$v_{\a_{j+1}}$, with $\lim_{j\to-\infty} v_{\a_{j}}=x$ and $\lim_{j\to\infty} v_{\a_{j}}=y$, the
limit being in the $R$-topology.  The length of the path is the sum of the resistances of the
constituent cells.  If $x$ (or $y$) is in $V_{\ast}$ we permit that the sequence begins with
infinite repetition of $x$ (respectively ends with repetition of $y$) which are considered
connected by the null cell of resistance zero, but otherwise the $v_{\a_{j}}$ are distinct.  Let
$S(x,y)$ denote the infimum of the lengths of paths from $x$ to $y$; it is easy to see from the
\frcs\ that there is a geodesic path that has length $S(x,y)$.

\BDF{df-loc-met}{We call the geodesic metric $S(x,y)$ the local resistance metric.}

\begin{lem}\label{lemma-compare-R-S}
For $x$ and $y$ in $\Omega$,
\begin{equation*}
    \frac{1}{2}S(x,y)
    \leq R(x,y)
    \leq S(x,y).
    \end{equation*}
\end{lem}

\Pr First consider the special case in which $x$ and $y$ are both in a loop-type cell $J_{\a}$, and
neither is contained in any smaller loop-type cell.  In this case none of the smaller loop-type
cells affects $R(x,y)$ or $S(x,y)$, so we may replace each such loop by its boundary vertex. The
result is to reduce the loop-type cell to a topological circle.  Deleting $x$ and $y$ from this
circle leaves two resistors, one with resistance $S(x,y)$ and the other with resistance at least
$S(x,y)$. The resistance $R(x,y)$ is the parallel sum of these, so satisfies $\frac{1}{2}S(x,y)\leq
R(x,y)\leq S(x,y)$.  This case also applies if both $x$ and $y$ are in $J_{(1)}\cup J_{(2)}$ and
neither is in any loop-type cell.

To complete the proof we show that the resistance from $x$ to $y$ decomposes as a series of
loop-type cells of the above form. Consider the (possibly empty) collection of loop-type cells that
contain $x$ but not $y$, and order them by inclusion, beginning with the largest.  Let
$v_{\a_{0}},v_{\a_{1}},v_{\a_{2}}\dotsc$ be the boundary vertices of these loops, and observe that
$v_{\a_{j}}\rightarrow x$.  Do the same for the loop-type cells containing $y$ but not $x$,
labeling the vertices $v_{\a_{-1}},v_{\a_{-2}},\dotsc$.  If $x$ is in $V_{*}$ then the sequence
will terminate  with an infinite repetition of $x$, and similarly for $y$.  Notice that deleting
any of the $v_{\a_{j}}$ disconnects $v_{\a_{j-1}}$ from $v_{\a_{j+1}}$.  This implies both that
that the effective resistances $R(v_{\a_{j}},v_{\a_{j+1}})$ sum in series to give the effective
resistance $R(x,y)$, and that the resistances $S(v_{\a_{j}},v_{\a_{j+1}})$ sum to $S(x,y)$.
However each of the configurations $v_{\a_{j}},v_{\a_{j+1}}$ is of the form of the special case
given above, so satisfies
\begin{equation*}
    \frac{1}{2}S(v_{\a_{j}},v_{\a_{j+1}})
    \leq R(v_{\a_{j}},v_{\a_{j+1}})
    \leq S(v_{\a_{j}},v_{\a_{j+1}}).
    \end{equation*}
Summing over $j$ then gives the desired inequality. \rP

By virtue of Theorem~\ref{Dirformdeterminedbytraces} and~\eqref{formonfiniteset} it is apparent
that we may describe a resistance form on $V_{*}$ in terms of the values $r_\a$.  The simple
structure of the graphs makes it easy to describe the choices of $\{r_{\a}\}_{\a}$ that give a
resistance form.

\begin{lem}\label{compatibilityofralphas}
Defining resistance forms on each $V_{n}$ by~\eqref{formonfiniteset} produces a sequence $\E_{n}$
that is compatible in the sense of Theorem~\ref{Dirformconstructedfromtraces} if and only if for
each arc-type cell $J_\a$,
\begin{equation*}
    r_{\a}=r_{\a1}+r_{\a2}.
    \end{equation*}
\end{lem}
\Pr%
Resistance forms satisfy the well-known Kirchoff laws from electrical network theory (see
\cite{Ki}, Section 2.1). If $J_{\a}$ is an arc-type cell then $J_{\alpha1}$ and $J_{\alpha2}$
connect $v_{\alpha1}$ and $v_{\alpha2}$ in series.  The resistance in $V_{|\a+1|}$ between
$v_{\alpha1}$ and $v_{\alpha2}$ neglecting $J\setminus J_{\a}$ is then $r_{\a1}+r_{\a2}$, so is
compatible with the resistance across $J_{\a}$ in $V_{|\a|}$ if and only if
$r_{\a}=r_{\a1}+r_{\a2}$.  In the alternative circumstance where $J_{\a}$ is a loop-type cell,
there is only one boundary vertex, so $r_\a$ is not defined and no constraint on $r_{\a1}$ and
$r_{\a2}$ is
necessary.%
\rP

According to Lemma~\ref{compatibilityofralphas} and Theorem~\ref{Dirformconstructedfromtraces}, one
may construct a resistance form on $V_{\ast}$ simply by choosing appropriate values $r_\a$.  It is
helpful to think of choosing these values inductively, so that at the $n$-th stage one has the
values $r_\a$ with $|\a|=n$.  In this case there are two types of operation involved in passing to
the $(n+1)$-th stage.  For arc-type cells $J_\a$ with $|\a|=n$ one chooses $r_{\a1}$ and $r_{\a2}$
so they sum to $r_\a$.  For loop-type cells $J_\a$ one chooses $r_{\a1}$ and $r_{\a2}$ freely.

This method provides a resistance form on $V_{*}$ and its $R$-completion $\Omega$, but our goal is
to describe Dirichlet forms on the fractal $J$.  In order that $\Omega=J$, or equivalently that the
topology from $\mathbb{C}$ coincides with the $R$-topology on $V_*$, we must further restrict the
values of $r_\a$. In the theorem below, $S-\diam{\bigl(O\bigr)}$ denotes the diameter of a set $O$
with respect to the local resistance metric $S(x,y)$.

\begin{thm}\label{thm-d}
The local regular resistance forms on $V_*$ for which $\Omega=J$ and the $R$-topology
is the same as the induced $\mathbb C$-topology are in one-to-one correspondence with the families
of positive numbers $r_{\a}$, one for each arc-type cell $J_{\a}$, that satisfy the conditions
\begin{gather}
    r_{\a}=r_{\a1}+r_{\a2} \label{alphacompateqn}\\
    \lim_{n\rightarrow\infty} \max_{\a\in\A_{n}} \Bigl(S-\diam{\bigl(J_{\a}\bigr)}\Bigr) = 0.
    \label{Sdiamsofcellsgotozero}
    \end{gather}

 A
sufficient but not necessary condition that implies~\eqref{Sdiamsofcellsgotozero}, and is often more
convenient, is
\begin{equation}\label{sufficientforSdiamsofcellsgotozero}
    \sum_{n} \max_{\a\in\A_{n}} r_{\a} <\infty.
    \end{equation}
\end{thm}

\Pr%
Lemma~\ref{compatibilityofralphas} shows that the condition~\eqref{alphacompateqn} on the $r_\a$ is
equivalent to compatibility of the sequence of resistance forms, which is necessary and sufficient
to obtain a resistance form on $\Omega$ by Theorems~\ref{Dirformdeterminedbytraces}
and~\ref{Dirformconstructedfromtraces}.

Recall that $V_*$ is $\mathbb{C}$-dense in the complete metric space $J$, so $J$ is the
$\mathbb{C}$-completion of $V_*$.  Similarly, $\Omega$ is by definition the $R$-completion of
$V_*$.   Then $\Omega=J$ and the $R$-topology is the same as the induced $\mathbb C$-topology if
and only if every $\mathbb{C}$-Cauchy sequence in $V_*$ is $R$-Cauchy and vice-versa.

Suppose there is an arc-type cell $J_\a$ with $r_{\a}=0$.  Then the sequence defined by
$x_{2j}=v_{\a1}$ and $x_{2j+1}=v_{\a2}$ is $R$-Cauchy but not $\mathbb{C}$-Cauchy.  Conversely
suppose there is a sequence that is $R$-Cauchy but not $\mathbb{C}$-Cauchy.  Compactness of $J$
allows us to select two distinct $\mathbb{C}$-limit points $x$ and $y$ and these will have
$R(x,y)=0$. Then $S(x,y)=0$ by Lemma~\ref{lemma-compare-R-S}, thus there is a non-trivial path
joining $x$ to $y$ such that $r_{\a_{j}}=0$ for all arc-type cells $J_{\a_{j}}$ on the path. It
follows that $R$-Cauchy sequences are $\mathbb{C}$-Cauchy if and only if $r_{\a}>0$ for all
arc-type cells $J_\a$.

An equivalence class of $\mathbb{C}$-Cauchy sequences is a point $x\in J$, and as noted in
Proposition~\ref{frcsandtopology}, $x$ is canonically associated to the nested sequence
$\bigl\{U_{n}(x)\bigr\}_{n\in\mathbb{N}}$, where $U_{n}(x)$ is the union of the $n$-cells
containing $x$.  Hence $\mathbb{C}$-Cauchy sequences are $R$-Cauchy if and only if for each $x$ the
resistance diameter of $U_{n}(x)$ goes to zero when $n\rightarrow\infty$. Clearly this is true if
\begin{equation}\label{unifcvgenceofRdiamstozero}
    \lim_{n\rightarrow\infty} \max_{\a\in\A_{n}} \Bigl(R-\diam{\bigl(J_{\a}\bigr)}\Bigr) = 0.
    \end{equation}
and conversely if~\eqref{unifcvgenceofRdiamstozero} fails then there is $\epsilon>0$ and for each
$n$ a cell $J_\a$ with $|\a|=n$ and $R$-diameter at least $\epsilon$, so compactness of $J$ gives a
$\mathbb{C}$-limit point $x$ at which the $R$-diameter of $U_{n}(x)$ is bounded below by $\epsilon$
independent of $n$. Applying Lemma~\ref{lemma-compare-R-S} we then see that $\mathbb{C}$-Cauchy
sequences are $R$-Cauchy if and only if~\eqref{Sdiamsofcellsgotozero} holds.

For any cell $J_\a$ and any $x\in J_\a$ there is a path from $v_{\a1}$ to $x$ which contains at
most one arc-type cell of each scale less than $|\a|$, so the
condition~\eqref{sufficientforSdiamsofcellsgotozero} implies~\eqref{Sdiamsofcellsgotozero}.  To see
this condition is not necessary we consider $\{r_{\a}\}$ as follows.  Fix a collection of loop-type
cells $J_\a$, one for each scale $|\a|\geq1$, with the property that if $J_\a$ is in this
collection then no loop-type ancestor of $J_\a$ is in the collection.  For example the cells with
addresses $(3),(13),(113),(1113),\dotsc$.  If $J_\a$ is in this collection set
$r_{\a1}=r_{\a2}=|\a|^{-1}$.  If $J_\a$ is a loop-type cell not in this collection set
$r_{\a1}=r_{\a2}=2^{-|\a|}$.  Also let $r_{\a1}=r_{\a2}=r_{a}/2$ if $J_\a$ is arc-type.  For these
values $r_\a$ we see that any local resistance path contains at most one arc-type cell from this
collection, so the $S$-diameter of a cell of scale $n\geq2$ is at most $(n-1)^{-1}+2^{-n-2}$
and~\eqref{Sdiamsofcellsgotozero} holds.  However $\max_{\a\in\A_{n}} r_{\a} = (n-1)^{-1}$ for each
$n\geq 2$, so~\eqref{sufficientforSdiamsofcellsgotozero} fails.\rP

\begin{cor}
Under the conditions of Theorem~\ref{thm-d}, all the functions in $\dom(\E)$ are continuous in the
topology from $\mathbb{C}$.
\end{cor}

\Pr%
It follows from (RF4) in Definition~\ref{d-resistance form}  that functions in $\dom(\E)$ are
$\frac{1}{2}$-H\"{o}lder continuous in the $R$-topology. \rP

The $n$-harmonic functions have a particularly simple form when written with respect to the local
resistance metric.

\BTHM{prop-loc-met}{An $n$-harmonic function is piecewise linear in the local resistance metric. }

\Pr The complement of $V_{n}$ is the finite union of cells $J_{\a}$ with $\a\in\A_{n}$.  If $f$ is
prescribed at $v_{\a1}$ and $v_{\a2}$ then its linear extension to $J_{\a}$ in the local resistance
metric has $f(v_{\a3})$ satisfying $r_{\a}f(v_{\a3})=r_{\a1}f(v_{\a2})+r_{\a2}f(v_{\a1})$.  However
the terms in the trace of $\E$ to $V_{n+1}$ that correspond to $J_{\a}$ are
\begin{equation*}
    r^{-1}_{\a1}\bigl( f(v_{\a1})-f(v_{\a3})\bigr)^{2} + r^{-1}_{\a2}\bigl( f(v_{\a2})-f(v_{\a3})\bigr)^{2}
    \end{equation*}
and it is clear this is minimized at precisely the given choice of $f(v_{\a3})$.  We have therefore
verified that an $n$-harmonic function extends from $V_{n}$ to $V_{n+1}$ linearly in the local
resistance metric, and the full result follows by induction. \rP

It is sometimes helpful to think of the local resistance metric as corresponding to a local
resistance measure $\nu$, which is defined as follows.

\begin{definition}
The  local
resistance measure of a compact set $E$ is given by
\begin{equation}\label{defnoflocresistmsr}
    \nu(E) = \inf \Bigl\{\sum_{j}r_{\a_{j}}: \bigcup_{j} J_{\a_{j}}\supset E\Bigr\}.
    \end{equation}
\end{definition}
One can easily see that $\nu$ has a unique extension to a positive, possibly infinite, Borel measure, and that
 $S(x,y)$ is the smallest measure of a path from $x$ to $y$.

This measure has a natural connection to the energy measures corresponding to functions in
$\dom(\E)$.  If $\E$ is local and $f\in\dom(\E)$, then the standard way to define the energy
measure $\nu_{f}$ is by the formula
$$ \int g\,d\nu_{f}=2\E (f,fg)-\E (f^2,g)$$
for any bounded quasi-continuous $g\in\dom(\E)$, see for instance \cite{FOT}. If $E$ is open, then
another   way to define $\nu_{f}(E)$ is to take the limit defining $\E$ from the resistance form as
in Theorem~\ref{Dirformconstructedfromtraces} and restrict to edges in $E$.  One may informally
think of the energy measure $\nu_{f}(E)$ of a set $E$ as being $\nu_{f}(E)=\E(f_{E})$, where
$f_{E}$ is equal to $u$ on $E$ and zero elsewhere, though this intuition is non-rigorous because
$f_{E}$ may fail to be in the domain of $\E$.   If $h_{m}$ is the piecewise harmonic function equal
to $u$ on $V_{m}$ then $\nu_{h_{m}}\rightarrow\nu_{f}$, and Theorem~\ref{prop-loc-met} ensures
$h_{m}$ has constant density $\frac{d\nu_{f}}{d\nu}$ on every sufficiently small cell. This
sequence of piecewise constant densities
is bounded by $\E(f)$ in $L^{1}(d\nu)$, and is a uniformly integrable submartingale.
The limit is the density of $\nu_{f}$ with respect to $\nu$, hence all energy measures are
absolutely continuous with respect to $\nu$.  If we let $\{\Omega_{j}\}$ be  the bounded Fatou
components of the polynomial $P(z)=z^{2}-1$  and note that $S(x,y)$ provides a local
parametrization of the topological circle $\partial\Omega_j$, then the above argument gives the
following description of the resistance form and the energy measures.

\begin{thm} \label{thm-der}
Under the conditions of Theorem~\ref{thm-d}, the domain $\dom(\E)$ of $\E$ consists of all
continuous functions such their restriction to each $\partial\Omega_j$ is absolutely continuous
with respect to the parametrization by the  local resistance metric, and the naturally defined
derivative $\frac{df}{dS}$ is square integrable with respect to $\nu$.  Moreover, each measure
$\nu_f$ is absolutely continuous with respect to $\nu$,
$$\frac{d\nu_{f}}{d\nu}=\left(\frac{df}{dS}\right)^2$$$\nu$-almost everywhere, and
$$\E(f,f)=\sum_{\Omega_j}\int_{\partial\Omega_j}\left(\frac{df}{dS}\right)^2d\nu.$$
\end{thm}

Note that the derivative $\frac{df}{dS}$ can be defined only up to orientation of the boundary
components $\partial\Omega_j$, but the densities and integrals in this theorem are independent of
this orientation. In general, $\nu$ is non-atomic, $\sigma$-finite, and for each $j$ we have
$$0<\nu(\partial\Omega_j)<\infty.$$
It should also be noted that $\nu$ plays only an auxiliary role in this theory, and is not
essential for the definitions of the energy or the \Lp.

Two specific choices of $\nu$ corresponding to resistance forms of the type described in
Theorem~\ref{thm-d} will be examined in more detail in Sections~\ref{specdecsection} and
\ref{sec-conf}. In both these cases $\nu$ is not finite, but $\sigma$-finite.

\section{\Lp s on the \bjs}\label{Lpsbjs}
As is usual in analysis on fractals, we use the Dirichlet form to define a weak Laplacian.  If
$\mu$ is a finite Borel measure on $J$, then the Laplacian with boundary behavior $B$ is defined by
\begin{equation}\label{defnofLaplacian}
    \E(f,g) = -\int_{J} (\Delta_{B} f) g \, d\mu \quad\text{for all $g\in\dom_{B}(\E)$}
    \end{equation}
where $\dom_{B}(\E)$ is the subspace of functions in $\dom(\E)$ satisfying the boundary condition
$B$.  In particular, if there is no boundary condition we have the Neumann Laplacian $\Delta_{N}$
and if the boundary condition is that $g\equiv0$ on $V_{0}$ we obtain the Dirichlet Laplacian
$\Delta_{D}$. We may then define a boundary operator $\partial_{n}^{B}$ such
that~\eqref{defnofLaplacian} can be extended to a general Gauss-Green formula.
\begin{equation}\label{gaussgreenfmla}
    \E(f,g) = -\int_{J} (\Delta_{B} f) g \, d\mu + \sum_{x\in V_{0}} g(x) \partial_{n}^{B}f(x)
    \quad\text{for all $g\in\dom(\E)$}.
    \end{equation}
Proofs of the preceding statements may be found in \cite{Ki4}.

The Laplacian may also be realized as a renormalized limit of Laplacians on the graphs $G_{n}$ by
using the method from \cite{Ki}. For $x\in V_{n}$ let $\psi_{x}^{n}$ denote the unique $n$-harmonic
function with $\psi_{x}^{n}(y)=\delta_{x,y}$ for $y\in V_{n}$, where $\delta_{x,y}$ is Kronecker's
delta. Since this function is $n$-harmonic, $\E(u,\psi_{x}^{n})=\E_{n}(u,\psi_{x}^{n})$ for all
$u\in\dom(\E)$. From this and \eqref{formonfiniteset} we see that if $x$ is in $V_{n-1}$ then
\begin{equation*}
    \E_{n}(u,\psi_{x}^{n})
    = \sum_{y\sim_{n} x} r^{-1}_{xy} \bigl( u(x)-u(y) \bigr),
    \end{equation*}
where $y\sim_{n}x$ indicates that $y$ and $x$ are endpoints of a common arc-type $n$-cell, and
$r_{xy}$ is the resistance of that cell.  We may view the expression on the right as giving the
value of a Laplacian on $G_{n}$ at the point $x$
\begin{equation}\label{laplacianonGn}
    \Delta_{n}^{r}u(x) = \sum_{y\sim_{n} x} r^{-1}_{xy} \bigl( u(x)-u(y) \bigr)
    \end{equation}
where the superscript $r$ in $\Delta_{n}^{r}$ indicates its dependence on the resistance form. By
the Gauss-Green formula~\ref{gaussgreenfmla},
\begin{equation*}
    \E_{n}(u,\psi_{x}^{n})= -\int (\Delta u) \psi_{x}^{n} \, d\mu
    \end{equation*}
so that
\begin{equation}\label{generalformofLaplacianaslimitofGnthings}
    \Bigl( \int \psi_{x}^{n}  \, d\mu \Bigr)^{-1} \Delta^{r}_{n}u(x)
    = \frac{-\int (\Delta u) \psi_{x}^{n} \, d\mu} {\int \psi_{x}^{n}  \, d\mu}
    \rightarrow -\Delta u(x)
    \end{equation}
as $n\rightarrow\infty$, which expresses $\Delta$ as a limit of the graph Laplacians
$\Delta_{n}^{r}$, renormalized by the measure $\mu$.

\section{Spectral decimation for a self-similar but not conformally invariant \Lp}\label{specdecsection}

The procedure in~\eqref{laplacianonGn} and~\eqref{generalformofLaplacianaslimitofGnthings} is
especially of interest when both the resistance form and the measure have a self-similar scaling
that permits us to express $\Delta^{r}_{n}$ in terms of the usual graph Laplacian
\begin{equation*}
    \Delta_{n}u(x) = \sum_{y\sim_{n} x} \bigl( u(x)-u(y) \bigr)
    \end{equation*}
and to simplify the expression for the measure.  Consider for example the simplest situation, in
which a resistance form is constructed on $J$ by setting $r_{\a}=2^{-|\a|}$, where $|\a|$ is the
length of the word $\a$ and using~\eqref{formonfiniteset}, and a Dirichlet form is obtained as in
Theorem~\ref{Dirformconstructedfromtraces}.  We take the measure $\mu_B$ to be the natural Bernoulli
one in which each $n$-cell has measure $(4\cdot3^{n-1})^{-1}$ for $n\geq1$.  In this
case~\eqref{laplacianonGn} simplifies to $\Delta^{r}_{n}u(x)=2^{n}\Delta_{n}u(x)$, and since $\int
\psi_{x}^{n}\,d\mu_B=2^{-1}3^{-n}$ we may reduce~\eqref{generalformofLaplacianaslimitofGnthings} to
\begin{equation}\label{defnofLaplacianaspointwiselimit}
    2\cdot 6^{n} \Delta_{n}u(x) \rightarrow - \Delta u(x).
    \end{equation}
The negative sign occurring on the right of~\eqref{defnofLaplacianaspointwiselimit} is a
consequence of the fact that  $\Delta_{n}$ is positive definite, whereas the
definition~\eqref{defnofLaplacian} gives a negative definite Laplacian.  The former is more
standard on graphs and the latter on fractals.

For the remainder of this section we study the particular \Lp\ defined in
\eqref{defnofLaplacianaspointwiselimit} using its graph approximations.  We begin by computing the
eigenstructure of the graph Laplacian on $G_{n}$ using the method of spectral decimation
(originally from~\cite{FukushimaShima,RammalToulouse,Shima}, though we
follow~\cite{Vibration,MalozemovTeplyaev}). The situation may be described as follows. The
transition matrix $M_{n}$ for a simple random walk on $G_{n}$ is an operator on the space of
functions on $V_{n}$. If we decompose this space into the direct sum of the functions on $V_{n-1}$
and its orthogonal complement, then $M_{n}$ has a corresponding block form
\begin{equation}\label{eq-M}
    M_{n}
    =\begin{pmatrix}
    A_{n}&B_{n}\\
    C_{n}&D_{n}
    \end{pmatrix}
\end{equation}
in which the matrix $A_{n}$ is a self-map of the space of functions on $V_{n-1}$. Define the Schur
complement $S$ to be  $A_{n}-B_{n} D_{n}^{-1} C_{n}$, and consider the Schur complement of the
matrix $M_{n}-z=M_{n}-zI$:
\begin{equation}\label{eq-Slambda}
S_{n}(z)=A_{n}-z-B_{n} (D_{n}-z)^{-1} C_{n}.
\end{equation}
If it is possible to  solve
\begin{equation}\label{eq-Slambda1}
S_{n}(z)=\phi_{n}(z)\Big(M_{n-1}-R_{n}(z)\Big),
\end{equation}
where $\phi_{n}(z)$ and $R_{n}(z)$ are scalar-valued rather than matrix-valued rational functions,
then we say that $M_{n}$ and $M_{n-1}$ are spectrally similar.  If we have a sequence $\{M_{n}\}$
in which each $M_{n}$ is a probabilistic graph \Lp\ on $G_n$ and $M_{n}$ is spectrally similar to
$M_{n-1}$, then it is possible to compute both the eigenvalues and eigenfunctions of the matrices
$M_{n}$ from $\phi_{n}(z)$ and $R_{n}(z)$.  Excluding the {\em exceptional set}, which  consists of the
eigenvalues of $D_{n}$ and the poles of $\phi_{n}(z)$, it may be shown that $z$ is an eigenvalue of
$M_{n}$ if and only if $R_{n}(z)$ is an eigenvalue of $M_{n-1}$, and the map $f\mapsto f- (D_{n}-z
)^{-1}C_{n}f$ takes the eigenspace of $M_{n-1}$ corresponding to $R_{n}(z)$ bijectively to the
eigenspace of $M_{n}$ corresponding to $z$ (\cite[Theorem 3.6]{MalozemovTeplyaev}).

Now consider a self-similar random walk on the graph $G_{n}$ in which the transition probability
from $v_{\a3}$ to $v_{\a1}$ in an arc-type cell is  a fixed number $p\in(0,1/2)$.  The transition
matrix for the cell $J_{\alpha}$ has the form
\begin{equation}\label{transitmatrixindecimationcase}
    M
    =\begin{pmatrix}
    1 & 0 & -1 \\
    0 & 1 & -1 \\
    -p & -p & 2p
    \end{pmatrix}.
\end{equation}
The results of \cite{Vibration} imply that the spectral decimation method is applicable to the
graph $G_{n}$. Moreover, self-similarity implies that both $\phi_{n}(z)$ and $R_{n}(z)$ are
independent of $n$  and may be calculated by examining a single cell $J_{\alpha}$.   From
\eqref{transitmatrixindecimationcase} we see that the eigenfunction extension map is
\begin{equation*}
    (D-z)^{-1}C
    =\begin{pmatrix}
    \frac{p}{2p-z } & \frac{p}{2p-z },
    \end{pmatrix}
\end{equation*}
meaning that the value at $v_{\a3}$ of a $\Delta_{|\a|+1}$ eigenfunction is $\frac{p}{2p-z}$ times
the sum of the values at $v_{\a1}$ and $v_{\a2}$. Since $M_{n-1}$ on a single cell is simply
\begin{equation*}
    M_{0}=\begin{pmatrix}
    1 & -1\\
    -1& 1
    \end{pmatrix},
    \end{equation*}
we find that
\begin{gather}
    \phi(z)=\frac{p}{2p -z },\quad\text{and} \notag\\
    R(z) = \frac{2p+1}p z-\frac1pz^2. \label{formulaforR}
    \end{gather}
The exceptional set is exactly the point $\{2p\}$. If we choose the initial \Lp\ on $G_0$ to be
\begin{equation*}
    \Delta_0=\begin{pmatrix}
    q & -q \\
    -q & q
\end{pmatrix}
\end{equation*}
for some $0<q<1$ then we can apply Proposition 4.1 and Theorem 4 of \cite{Vibration} to compute
both the multiplicities and the eigenprojectors.

\def\mult#1#2{\textup{mult}_{#1}{(#2)}}

\BTHM{thm-Gn}{The eigenvalues of the Laplacian $\Delta_{n}$ on $G_{n}$ are given by
\begin{gather*}
    \sigma(\Delta_0) = \{0,2q\},\\
    \sigma(\Delta_n)=   \left( \bigcup_{m=0}^{n-1} R^{-m}\{2p\}\right) \bigcup \left(R^{-n}\{0,2q\}\right).
    \end{gather*}
Moreover, if $z\in R^{-n}\{0,2q\}$ then $\mult nz = 1$ and the corresponding eigenfunctions have
support equal to $J$; if $z\in R^{-m}\{2p\}$ then $\mult n{z} = 2 \cdot 3^{n-m-1}$ and the
corresponding eigenfunctions vanish on $V_{n-m-1}$.}

\Pr For $z\in R^{-n}\{0,2q\}$ the result follows from Proposition 4.1(i) and Theorem 4(i) of
\cite{Vibration}. For $z\in R^{-m}\{2p\}$ the result follows from Proposition 4.1(iii) and Theorem
4(iii) of \cite{Vibration}. In particular,
\begin{equation*}
    \mult n{2p}=4\cdot3^{n-1}-|V_{n-1}|+\mult{n-1}{R(2p)},
    \end{equation*}
where $R(2p)=2$, which is not in the spectrum of any $\Delta_k$, and $|V_k|=2\cdot3^{k}$. \rP

\BCOR{cor}{The normalized limiting distribution of eigenvalues (also called the integrated density
of states) is a pure point measure $\kappa$ with atoms at each point of the set
$$  \bigcup_{m=0}^{\infty} R^{-m}\{2p\}, $$
Moreover, if $ z\in R^{-m}\{2p\}$ then $\kappa (\{z\})=2\cdot3^{-m-1}$. There is one atom in each
gap of the Julia set of $R$. }

A special case occurs if we make the convention that every edge can be traveled in both directions
with equal probability, in which case each of the $G_{n}$ is a regular graph of degree 4. This
simple random walk has $p=q=\frac14$ from which $R(z)=6z-4z^2$. Since our graphs have $2\cdot3^{n}$
vertices, we conclude that in this case the spectral dimension of the corresponding infinite graphs
is
\begin{equation*}
    d_s=\frac{2\log3}{\log6}.
    \end{equation*}
One may also consider weighted Laplacians on the infinite graphs by varying the parameter $p$.

We saw at the end of Section~\ref{Lpsbjs} that the Laplacian $\Delta$ on the fractal $J$ may be
obtained as a limit of graph Laplacians $\Delta_{n}$, provided that both the Dirichlet form and the
measure have self-similar scaling.  Under these circumstances, the spectral decimation method gives
a natural algorithm for constructing eigenfunctions of the Laplacian on the fractal. This method
was first developed for the Sierpinski Gasket fractal \cite{RammalToulouse,Shima,FukushimaShima}.

We illustrate this method for the special self-similar case where the resistance form on $J$
satisfies~\eqref{formonfiniteset} with $$r_{\a}=2^{-|\a|},$$ where $|\a|$ is the length of the word
$\a$, and the Dirichlet form is obtained using Theorem~\ref{Dirformconstructedfromtraces}.  In this
case
\begin{equation*}
    \E_{n}(u,\psi_{x}^{n})
    = \sum_{y\sim_{n} x} r^{-1}_{xy} \bigl( u(x)-u(y) \bigr)
    = 2^{n} \sum_{y\sim_{n} x}  \bigl( u(x)-u(y) \bigr)
    =4\cdot2^{n} \Delta_{n}u(x)
    \end{equation*}
where $\Delta_{n}$ is the graph Laplacian on $G_{n}$ with equal weight $\frac14$ on each edge. This
is equivalent to setting $p=q=\frac14$. Correcting for the extra
factor of $\frac14$ in the graph Laplacian we find from~\eqref{defnofLaplacianaspointwiselimit}
\begin{equation}\label{Laplacianaspointwiselimitinspecdeccase}
    8\cdot 6^{n} \Delta_{n}u(x) \rightarrow - \Delta u(x).
    \end{equation}
Here we take that the measure $\mu$ in~\eqref{defnofLaplacian} is the
 the natural Bernoulli measure $\mu_B$
for which each $n$-cell has $\mu_B$-measure equal to $(4\cdot3^{n-1})^{-1}$ for $n\geq1$.

Now suppose that $\{u_{n}\}$ is a sequence of eigenfunctions of  $\Delta_{n}$ with eigenvalues
$\lambda_{n}$, and the property that $u_{n}=u_{m}$ on $V_{m}$ for $m\leq n$. Further assume that
$6^{n}\lambda_{n}$ converges and that the function $u$ defined on $V_{\ast}$ by $u(x)=u_{n}(x)$ for
$x\in V_{n}$ is uniformly continuous, and thus can be extended continuously to $J$.
Then~\eqref{Laplacianaspointwiselimitinspecdeccase} implies that $u$ is a Laplacian eigenfunction
on $J$ with eigenvalue $\lambda=-8\lim 6^{n}\lambda_{n}$. From the formula~\eqref{formulaforR} for
$R$ we have
\begin{equation*}
    \lambda_{n} = \frac{3+\epsilon_{n}\sqrt{9-4\lambda_{n-1}} }{4}
    \end{equation*}
where $\epsilon_{n}$ is one of $\pm1$ for each $n$. If only finitely many $\epsilon_{n}$ equal $+1$
then $6^{n}\lambda_{n}$ converges and it is easily verified that $u$ is uniformly continuous on
$V_{\ast}$, so this method constructs a large number of eigenfunctions.  It is actually the case
that it constructs all eigenfunctions, though we will only show this for the Dirichlet Laplacian.

The Dirichlet eigenfunctions corresponding to $R^{-m}(2p)=R^{-m}(\frac{1}{2})$ produce Dirichlet
eigenfunctions on $J$.  Via an argument from \cite{FukushimaShima}, this provides a precise
description of the spectrum of the Dirichlet Laplacian $\Delta_{D}$.  Let $\psi(x)=
\frac{3-\sqrt{9-4x} }{4}$ and
\begin{equation*}
    \Psi(x)
    = \lim_{n\to\infty} 6^{n} \psi^{n}(x)
    \end{equation*}
in which the limit is well-defined on a neighborhood of zero by the Koenig's linearization theorem
(see \cite{Milnor}).  Note that $\Psi(0)=0$ and $\Psi'(0)=1$, so that $\Psi$ is also invertible on
a neighborhood of $0$.  In the above construction of Dirichlet eigenvalues we asked that all but
finitely many of the inverse branches of $R$ be exactly $\psi$, so that for any such
$\lambda=-8\lim 6^{n}\lambda_{n}$ there is $n_{0}$ such that $\lambda_{n+1}=\psi(\lambda_{n})$ for
all $n\geq n_{0}$.  It follows that
\begin{equation*}
    \lambda
    =-8\cdot \lim_{n\rightarrow\infty} 6^{n_{0}} 6^{n-n_{0}} \psi^{n-n_{0}}(\lambda_{n_{0}})
    =-8\cdot 6^{n_{0}} \Psi(\lambda_{n_{0}})
    \end{equation*}
where $\lambda_{n_{0}}=R^{-m}\bigl(\frac{1}{2}\bigr)$ for some $0\leq m\leq n_{0}$.

\begin{thm}
The spectrum of $\Delta_{D}$ on $J$ consists of isolated eigenvalues 
$$\lambda=-8\cdot6^{n_{0}}\Psi\bigl(R^{-m}(\frac{1}{2})\bigr)$$ with multiplicity $2\cdot3^{n_{0}-m-1}$, for
each $n_{0}\geq1$ and $0\leq m\leq n_{0}$.  The corresponding eigenfunctions are those obtained
from the eigenfunctions in Theorem~\ref{thm-Gn} by spectral decimation.
\end{thm}

\begin{proof}
Kigami~\cite{Ki4} proves that there is a Green's operator with a kernel that is uniformly Lipschitz
in the resistance metric, hence the resolvent of the Laplacian is compact and the spectrum of the
\Lp\ is discrete (pure point with isolated eigenvalues of finite multiplicity accumulating to infinity).  Since $\Delta_{D}$ is negative definite, the spectrum consists of a
decreasing sequence $\lambda_{j}$ of negative real eigenvalues that accumulate only at $-\infty$.

We have seen that the spectral decimation construction produces some Dirichlet eigenvalues and
their eigenfunctions.  The standard way to determine that all points in the spectrum arise in this
manner is a counting argument due to Fukushima and Shima~\cite{FukushimaShima}.  As the argument
holds essentially without alteration, we only sketch the details.

Expanding the Green's kernel $g(x,y)$ of the Laplacian as an $L^{2}$-series in the eigenfunctions,
we find that
\begin{equation*}
    - \int g(x,x) d\mu_B(x) = \sum_{i} \frac{1}{\lambda_{i}}
    \end{equation*}
where the sum is over the eigenvalues of $\Delta_{D}$, each repeated according to its multiplicity.
Similarly, if we let $g_{m}$ be the Green's kernel for $-8\cdot6^{m}\Delta_{m}$ and let $\mu_{m}$
be the measure with equal mass at each point of $V_{m}$, then
\begin{equation*}
    - \int g_{m}(x,x) d\mu_{m}(x) = \sum_{j} \frac{1}{\kappa^{(m)}_{j}}
    \end{equation*}
where the sum is over its eigenvalues.  However $g$ is continuous and equal to $g_{m}$ on $V_{m}$,
and the measures $\mu_{m}$ converge $\text{weak}^{*}$ to $\mu_B$, so as $m\rightarrow\infty$ the sum
of all $\frac{1}{\kappa^{(m)}_{j}}$ converges to the sum of $\frac{1}{\lambda_{i}}$.

Now each $\kappa^{(m)}_{j}$ is $-8\cdot6^{m}\lambda_{j}^{(m)}$, where $\lambda_{j}^{(m)}\in
R^{-m}(\frac{1}{2})$, and any sequence $-8\cdot6^{m}\lambda_{j}^{(m)}$ satisfying the conditions of
the spectral decimation algorithm converges to some eigenvalue $\lambda_{i}$ of $\Delta_{D}$.  With
a little care it is possible to show that $\sum_{j} \frac{1}{\kappa^{(m)}_{j}}$ converges to the
sum $\sum_{k} \frac{1}{\lambda_{i_{k}}}$, over those eigenvalues that arise from the spectral
decimation.  Since $\sum_{j} \frac{1}{\kappa^{(m)}_{j}}$ also converges to $\frac{1}{\lambda_{i}}$,
we conclude that the spectral decimation produces all eigenvalues.
\end{proof}

It is worth noting that eigenfunctions also have a self-similar scaling property. Specifically, let
$f_{\a}$ denote the natural map from $J_{(1)}$ to $J_{\a}$ if $J_{\a}$ is an arc-type cell, and
from $J_{(3)}$ to $J_{\a}$ if $J_{\a}$ is a loop-type cell.  This natural map is defined in the
obvious way on the boundary points and then inductively extended to map $V_{n}\cap J_{(1)}$ to
$V_{n+|\a|-1}\cap J_{\a}$ (respectively $V_{n}\cap J_{(3)}$ to $V_{n+|\a|-1}\cap J_{\a}$) for each
$n$, whereupon it is defined on the entire cell by continuity.  By the definition of the Dirichlet
form, this composition scales energy by $2^{1-|\a|}$, and
by~\eqref{defnofLaplacianaspointwiselimit} it scales the Laplacian by $6^{1-|\a|}$. More precisely,
if $u$ is such that $(\Delta-\lambda)u=0$ then $\Delta(u\circ f_{\a})=6^{1-|\a|}(\Delta u)\circ
f_{\a}$, so $\Delta(u\circ f_{\a})$ is a Laplacian eigenfunction with eigenvalue
$6^{1-|\a|}\lambda$.

The scaling property provides a very simple description of the Dirichlet eigenfunctions.  Suppose
$u$ is a Dirichlet eigenfunction obtained as the limit of $u_{n}$ according to the spectral
decimation, and let $m$ be the scale with $\lambda_{m}=\frac{1}{2}$.   Then $u_{m}$ vanishes on
$V_{m-1}$, so if $|\a|=m$ then $u_{m}\circ f_{\a}$ is a Dirichlet eigenfunction on $J_{(1)}$ (or
$J_{(3)}$) with eigenvalue $6^{1-m}$.  There is a one dimensional space of such functions (note
that whether the function is on $J_{(1)}$ or $J_{(3)}$ is immaterial because it vanishes on the
boundary), spanned by the Dirichlet eigenfunction on $J_{(1)}$ with value $1$ at $v_{13}$.  It
follows that the Dirichlet eigenfunctions are all built by gluing together multiples of this
function on individual cells of a fixed scale $m$, subject only to the condition that the values on
$V_{m}$ give a graph eigenfunction with eigenvalue $\frac{1}{2}$.

\section{Conformally invariant resistance form and \Lp}\label{sec-conf}

In this section we decompose $J$ as a union of a left and right piece $J=J_L\cup J_R$, where
\begin{gather*}
    J_L=J\cap\{z:Re(z)\leqslant \tfrac{1-\sqrt5}2\}=J_{(3)}\\
    J_R=J\cap\{z:Re(z)\geqslant \tfrac{1-\sqrt5}2\}=J_{(1)}\cup J_{(2)}\cup J_{(4)}.
    \end{gather*}
The sets meet at $a=\tfrac{1-\sqrt5}2$, which is the fixed point of $P(z)=z^2-1$.
The polynomial $P(z)$ maps $J_L$ onto $J_R$ by an one-to-one mapping, and the piece $J_{(4)}\subset
J_R$ onto $J_R$ by a one-to-one mapping. It also maps the central part $J_{(1)}\cup J_{(2)}$ of
$J_R$ onto $J_L$ by a two-to-one mapping. Therefore the following directed graph
\begin{equation*}
    \xymatrix{ J_L \ar@/_/[rr] &&  \ar@/_/[ll] \ar@/_/@<-.7ex>[ll] \ar@(ur,dr)[] J_R}
    \end{equation*}
corresponds to the action of $P(z)$, and defines a graph directed cell structure on $J$.  Note that
$V_{\ast}=\cup_{m} P^{-m}\{a\}$ and that the preimages of arc-type cells under $P$ are also
arc-type cells, while the preimages of loop-type cells are loop-type cells except in the case of
$J_{(3)}=J_{L}$ for which the preimages are $J_{(1)}$ and $J_{(2)}$. This construction is related
to group-theoretic results about these graphs and \cite[and references
therein]{Nekrashevych,NekandTep}, and in particular to the  substitution scheme in
Figures~\ref{d-c.i.g.} and~\ref{f-sub}, in which the labeling of components is $J_L=A$ and $J_R=B$.

As usual, we are interested in Dirichlet forms and measures that have a self-similar scaling under
natural maps of the fractal.  In this case, the mapping properties described above show that if
$\E$ is a Dirichlet form on $J$ then we may define Dirichlet forms $\E^{i}$ on the cells $J_{(i)}$,
$i=1,2,3,4$ by setting
\begin{equation*}
    \E^{i} (u) = \E(u\circ P) \quad\text{for $u$ on $J_{(i)}$ with $u\circ P\in\dom(\E)$}
    \end{equation*}
where $u\circ P$ is taken to be zero off $P(J_{(i)})$ in each case.  The form $\E$ is then
self-similar under the action of $P$ if for $u\in\dom(\E)$
\begin{equation}\label{selfsimilarunderP}
    \E(u) = \rho \sum_{i} \E^{i}\Bigl(u|_{J_{(i)}}\Bigr)
    \end{equation}
for some $\rho$.

\begin{thm}\label{invarresistformexists}
Among the resistance forms identified in Theorem~\ref{thm-d} there is one that has a self-similar
scaling under the action of $P(z)$ and is symmetric under complex conjugation.  It is unique up to
a scalar multiple, and has scaling factor $\rho=\sqrt{2}$.
\end{thm}

\begin{proof}
By Theorems~\ref{Dirformdeterminedbytraces} and~\ref{Dirformconstructedfromtraces} a necessary and
sufficient condition for~\eqref{selfsimilarunderP} to be true is that the trace of both $\E$ and
$\sum_{i}\E^{i}$ to $V_{m}$ are equal for each $m$.  The trace of $\E$ to $V_{m}$ is a resistance
form
\begin{equation*}
    \E_{m}(u) = \sum_{\a\in \A_{m}} r^{-1}_{\a} \bigl( u(v_{\a1})-u(v_{\a2})\bigr)^{2}
    \end{equation*}
as in~\eqref{formonfiniteset}. The trace of $\sum_{i}\E^{i}$ to $V_{m}$ is found by minimizing the
energy when values on $V_{m}$ are fixed, and each $\E^{i}$ may be minimized separately.  Thus for
each $i$ the restriction of $u$ to $J_{(i)}$ has the property that $u\circ P$ is energy minimizing
on $P(J_{(i)})$.  The result is therefore a resistance form in which the resistance across an
arc-type cell $J_{\a}$ is equal to the resistance of the form $\E$ across $P(J_{(\a)})$.

We conclude that~\eqref{selfsimilarunderP} is true if and only if $\E$ is the limit of resistance
forms $\E_{m}$ with $r_{P(J_{(\a)})}=\rho r_{J_{(\a)}}$.  There is only one value of $\rho$ for
which this can be satisfied.  To see this, note that $P^{2}$ maps both $J_{(11)}$ and $J_{(22)}$ to
$J_{(1)}$, and both $J_{(12)}$ and $J_{(21)}$ to $J_{(2)}$, so $r_{(11)}=r_{(22)}=\rho^{-2}r_{(1)}$
and $r_{(12)}=r_{(21)}=\rho^{-2}r_{(2)}$.  However $r_{(11)}+r_{(12)}=r_{(1)}$ and
$r_{(21)}+r_{(22)}=r_{(2)}$, so $r_{(1)}=r_{(2)}$ and $\rho^{2}=2$.  Also $r_{(i1)}$ and $r_{(i2)}$
are equal to $\frac{1}{2}r_{(1)}$ for $i=1,2$.

Observe that for any arc-type cell $J_{(\a)}$ there is a unique $m=m(\a)$ so
$P^{m}(J_{(\a)})=J_{(i)}$ for one of $i=1,2$.  According to our reasoning thus far, we must have
$r_{(1)}=r_{(2)}$ and  $r_{\a}=2^{-m(\a)/2}r_{(1)}$.  It remains to be seen that these resistances
satisfy the conditions of Theorem~\ref{thm-d}.  The condition $\lim_{|\a|\rightarrow\infty}
r_{\a}=0$ is immediate, and one can easily verify \eqref{Sdiamsofcellsgotozero}, in particular by the computation in Theorem~\ref{thm-J-res}. For any $\a$ we have  $P^{m(\a)}\bigl(J_{(\a1)}\bigr)$ and
$P^{m(\a)}\bigl(J_{(\a2)}\bigr)$ are $J_{(i1)}$ and $J_{(i2)}$ for one of $i=1,2$, so the second
condition $r_{\a1}+r_{\a2}=r_{\a}$ is equivalent to $r_{11}+r_{12}=r_{1}$, and the latter has
already been established.
\end{proof}

Recall that we have a local resistance metric  and a measure $d\nu$ (see
Definition~\ref{df-loc-met}, \eqref{defnoflocresistmsr} and Theorem~\ref{thm-der}) corresponding to
a sequence of resistance forms.  For the resistance form from
Theorem~\ref{invarresistformexists} the measure $d\nu$ is related to the harmonic measure on Fatou
components. In this case the measure is infinite.  It is well known that in our situation each Fatou component $\Omega\subset\mathbb{C}$
is a topological disc with locally-connected boundary, so the Riemann map $\sigma$ from the unit
disc to $\Omega$ has a continuous extension to the unit circle.  The harmonic measure from the
point $x\in\Omega$ is the image on $\partial\Omega$ of the Lebesgue measure on the circle under the
Riemann map with $0\mapsto x$.  It is the same as the exit probability measure for a Brownian
motion in $\Omega$ started at $x$, and is the representing measure for the linear functional on
$C(\partial\Omega)$ that takes $f$ to the value of the harmonic extension of $f$ at $x$.

\begin{thm}\label{thm-J-res}
Let $d\nu$ be the measure corresponding to the unique resistance form from
Theorem~\ref{invarresistformexists} with normalization $r_{(1)}=\frac12$. Then $\nu$ is an infinite measure that has
the self-similar scaling $\nu(P(E))=\sqrt2\nu(E)$ for any set $E$ on which $P$ is injective.

Let $\Omega=\Omega_{0}$ be the Fatou component of $P$ that contains the critical point $0$, and for
each other bounded component $\Omega_{j}$ of the Fatou set of $P$ let $m_{j}$ be the unique number
such that $P^{m_{j}}$ maps $\Omega_{j}$ bijectively to $\Omega$.  Let $d\nu_{j}$ be the harmonic
measure on $\partial\Omega_{j}$ from the point $P^{-m_{j}}(0)\in\Omega_{j}$.  Then
\begin{equation*}
    d\nu = \sum_{j=0}^{\infty} 2^{-m_{j}/2} d\nu_{j}
    \end{equation*}
\end{thm}

\begin{proof}
Let $\sigma$ be a Riemann map from the unit disc to $\Omega$ with $\sigma(0)=0$. Since $P^{2}$ is a
two-to-one map of $\Omega$ onto itself, $\sigma^{-1}\circ P^{2}\circ\sigma$ is a two-to-one map of
the unit disc onto itself. A version of the Schwartz lemma (for example, that in \cite{MR1691000})
implies that $\sigma^{-1}\circ P^{2}\circ\sigma=cz^{2}$, where $c$ is a constant with $|c|=1$.
Moreover, there is a unique $\sigma$ such that $\sigma(1)=a$, and so $\sigma^{-1}\circ
P^{2}\circ\sigma=z^{2}$. Since $\Omega$ is locally connected the Riemann map extends to the
boundary. Pulling back the measure $d\nu$ to the circle via $\sigma$ gives a Borel measure that
scales by $2$ under $z\mapsto z^{2}$. Consider the set of $2^m$ preimages of $a$ under the
composition power $P^{2m}$ that lie in $\partial\Omega$. These preimages divide
$\nu\big|_{\partial\Omega}$ into $2^m$ equal parts. The preimages of these $2^m$ points under
$\sigma$ are binary rational points on the unit circle that divide the Lebesgue measure into $2^m$
equal parts. It follows that $d\nu$ is a multiple of the harmonic measure for the point
$\sigma(0)=0$, and since they both have measure $1$ they are equal. A similar alternative
construction is to consider an ``internal ray'' which is the intersection of the negative real
half-line with $\Omega$, and its preimages. Then the harmonic measure can be determined in the
usual way by computing angles between these rays.

The bounded Fatou components of $P$ are $\Omega$ and the topological discs enclosed by loop-type
cells.  The argument we have just given applies to any such component $\Omega_{j}$, except that the
Riemann map is $P^{-m_{j}}\circ\sigma$, so $d\nu\bigl|_{\Omega_{j}}$ is a multiple of the harmonic
measure $d\nu_{j}$ from the point $P^{-m_{j}}(0)$.  The result then follows from the proof of
Theorem~\ref{invarresistformexists}, where it is determined that
$\nu(\partial\Omega_{j})=2^{-m_{j}/2}$.
\end{proof}

It is natural to compare this to other measures on $J$.  We saw in Theorem~\ref{thm-der} that the
energy measures are absolutely continuous to $d\nu$. Another standard measure to consider is the
unique  balanced invariant probability measure of $P$, denoted $\mu_P$. It can be obtained, for
instance, as the weak limit of the sequence of probability measures $\mu_{m}$, where each $\mu_{m}$
is $2^{-m}$ times the counting measure on the $2^{m}$ preimages of $a$. An alternative construction
of $\mu_P$ defines it as the harmonic measure from infinity, which also can be determined in the
usual way by computing angles between the external rays. This measure is a Bernoulli-type measure
that has the self-similar scaling $\mu_P(P(E))=2\mu_P(E)$ for any set $E$ on which $P$ is injective
(or any set if we incorporate multiplicity). The measures $\mu_P$ and $\nu$ are singular, as may be
verified directly by comparing $\mu_{P}$ to $\nu$. Indeed, $\mu_{P}$ has measure $2^{-m}$ on the
preimages $P^{-m}(J_{(i)})$, $i=1,2$ whereas $\nu(P^{-m}(J_{(i)}))=2^{-m/2}$.

Let  $\Delta_P$ be the \Lp\ corresponding to the unique conformally invariant $\E$ and the balanced
invariant measure $\mu_P$. Because of \Thm{invarresistformexists}, $\Delta_P$ is (up to a constant
multiple) the only \Lp\ that has self-similar scaling under the action of $P(z)=z^{2}-1$, and its
scaling factor is $2\sqrt2$.

\BTHM{thm-J-Lp}{The spectral dimension of $\Delta_P$ is equal to $\dfrac43$. }

\Pr Since  $J$ has graph-directed fractal structure,
 the method of \cite{HN,KL1} is
applicable.  This reduces the spectral dimension computation to finding $s$ such that the spectral
radius of the matrix
 $$
 \left(2\sqrt2\right)^{-s}
 \left[\begin{array}{cc}
0&2\\
1&1
\end{array}\right]
$$
is equal to one. Thus $s=\frac23$ and $d_s=2s=\frac43$.\rP

\subsection*{Acknowledgements}
The authors would like to thank John Hubbard, Jun Kigami, Volodymyr Nekrashevych and Robert Strichartz for several important
suggestions.  Some of the work on this paper was done during the Analysis on Graphs and its
Applications program at the Isaac Newton Institute for Mathematical Sciences, with funding provided
by the Issac Newton Institute, the London Mathematical Society and the National Science Foundation;
the authors are grateful to these institutions for their support.

\end{document}